\documentclass[]{article}

\usepackage[a4paper,
            top=1.25in,
            bottom=1.25in,
            left=1.25in,
            right=1.25in]{geometry}

\usepackage{siunitx}
\usepackage{bm}
\usepackage{amsmath,amsfonts,amssymb,amsthm}
\usepackage{graphicx}
\usepackage{epstopdf}
\usepackage{algorithmic}
\usepackage{booktabs}
\usepackage{multirow}

\usepackage{xcolor}
\definecolor{darkblue}{rgb}{0.0, 0.0, 0.5}
\definecolor{darkgreen}{rgb}{0.0, 0.4, 0.0}
\definecolor{darkred}{rgb}{0.5, 0.0, 0.0}
\usepackage[
    colorlinks=true,
    linkcolor=darkblue,       
    citecolor=darkgreen,      
    urlcolor=darkred,         
    filecolor=darkblue,       
    pdfborder={0 0 0}         
]{hyperref}

\usepackage[capitalize,nameinlink]{cleveref}
\crefname{section}{section}{sections}
\crefname{subsection}{subsection}{subsections}
\Crefname{section}{Section}{Sections}
\Crefname{subsection}{Subsection}{Subsections}
\Crefname{figure}{Figure}{Figures}
\crefformat{equation}{\textup{#2(#1)#3}}
\crefrangeformat{equation}{\textup{#3(#1)#4--#5(#2)#6}}
\crefmultiformat{equation}{\textup{#2(#1)#3}}{ and \textup{#2(#1)#3}}
{, \textup{#2(#1)#3}}{, and \textup{#2(#1)#3}}
\crefrangemultiformat{equation}{\textup{#3(#1)#4--#5(#2)#6}}%
{ and \textup{#3(#1)#4--#5(#2)#6}}{, \textup{#3(#1)#4--#5(#2)#6}}{, and \textup{#3(#1)#4--#5(#2)#6}}
\Crefformat{equation}{#2Equation~\textup{(#1)}#3}
\Crefrangeformat{equation}{Equations~\textup{#3(#1)#4--#5(#2)#6}}
\Crefmultiformat{equation}{Equations~\textup{#2(#1)#3}}{ and \textup{#2(#1)#3}}
{, \textup{#2(#1)#3}}{, and \textup{#2(#1)#3}}
\Crefrangemultiformat{equation}{Equations~\textup{#3(#1)#4--#5(#2)#6}}%
{ and \textup{#3(#1)#4--#5(#2)#6}}{, \textup{#3(#1)#4--#5(#2)#6}}{, and \textup{#3(#1)#4--#5(#2)#6}}

\usepackage{subcaption}
\usepackage{xspace}

\newcommand{\podrbf}{\ensuremath{\operatorname{POD\text{-}RBF}}\xspace}
\newcommand{\quadls}{\ensuremath{\operatorname{Quad-LS}}\xspace}
\newcommand{\quadnet}{\ensuremath{\operatorname{QuadNet}}\xspace}
\newcommand{\quadnetmu}{\ensuremath{\operatorname{QuadNet}\text{-}\mu}\xspace}
\newcommand{\RA}[1]{#1}
\newcommand{\RB}[1]{#1}
\newcommand{\RALL}[1]{#1}

\newcommand{\email}[1]{\protect\href{mailto:#1}{#1}}

\theoremstyle{remark}
\newtheorem{remark}{Remark}

\newcommand{\ntrain}{\ensuremath{N^{\text{train}}_{\mu}}}

\newcommand{\ndof}{\ensuremath{N_{dof}}}

\title{\textbf{Machine Learning-based quadratic closures for non-intrusive Reduced Order Models}

\footnotetext{
    \textbf{Funding:} This work was partially funded by PRIN ``Numerical Analysis for Full and Reduced Order Methods for Partial Differential Equations" (NA-FROM-PDEs) project and INdAM-GNCS: Istituto Nazionale di Alta Matematica –– Gruppo Nazionale di Calcolo Scientifico. This work has been conducted within the research activities of the consortium iNEST (Interconnected North-East Innovation Ecosystem), Piano Nazionale di Ripresa e Resilienza (PNRR) – Missione 4 Componente 2, Investimento 1.5 – D.D. 1058 23/06/2022, ECS00000043, supported by the European Union's NextGenerationEU program.}
    }

\author{Gabriele Codega\thanks{International School for Advanced Studies, SISSA, Trieste, Italy (\email{gcodega@sissa.it}).}
\and Anna Ivagnes\thanks{International School for Advanced Studies, SISSA, Trieste, Italy
  (\email{aivagnes@sissa.it}).}
\and Nicola Demo\thanks{International School for Advanced Studies, SISSA, Trieste, Italy; FAST Computing, Trieste, Italy (\email{ndemo@sissa.it}).}
\and Gianluigi Rozza\thanks{International School for Advanced Studies, SISSA, Trieste, Italy (\email{grozza@sissa.it).}}
}
\date{\vspace{-5ex}}

\begin{document}

\maketitle
\begin{abstract}
\noindent
In the present work, we introduce a data-driven approach to enhance the
accuracy of non-intrusive Reduced Order Models (ROMs). In particular,
we focus on ROMs built using Proper Orthogonal Decomposition (POD) in
an under-resolved and marginally-resolved regime, i.e. when the number of
modes employed is not enough to capture the system dynamics. We propose a method to re-introduce the contribution of neglected modes through a quadratic correction term, given by the action of a quadratic operator on the POD coefficients. Differently from the state-of-the-art methodologies, where the operator is learned via least-squares optimisation \cite{geelen2023quadratic-correction,barnett2022quadratic}, we propose to parametrise the operator by a Multi-Input Operators Network (MIONet). This way, we are able to build models with higher generalisation capabilities, where the operator itself is continuous in space -- thus agnostic of the domain discretisation -- and parameter-dependent. We test our model on two standard benchmarks in fluid dynamics and show that the correction term improves the accuracy of standard POD-based ROMs.
\end{abstract}

\section{Introduction}
\label{sec:intro}

The increasing availability of computational resources in the last decades has lead to the development of computational fluid dynamics (CFD).
Large methodological improvements have been made to accurately capture the flow fields, for example with large-eddy simulations (LES)~\cite{piomelli1999large, sagaut2006large}, which only resolve the largest scales, or Reynolds-averaged Navier--Stokes simulations (RANS) \cite{reynolds1895iv}, where only the time-averaged flow is resolved and the oscillations are modelled through appropriate turbulence modelling.
Despite the large availability of CPU resources, these types of simulations are still prohibitive in terms of computational time, especially in real-world and industrial applications. Indeed, many engineering and scientific applications, such as design optimization, decision making, control and uncertainty quantification, require fast predictions.

A large set of tools named reduced order models (ROMs)~\cite{degruyter1, degruyter2, degruyter3, quarteroni2014reduced, quarteroni2015reduced, rozza2022advanced} has been developed in the past years to accelerate simulations, while maintaining a good accuracy. ROMs are usually based on an \emph{offline-online} procedure. The \emph{offline} stage typically consists in collecting data from pre-computed expensive numerical simulations (like LES or RANS), while the \emph{online} stage provides an accurate prediction for unseen configurations in much less computational time, without re-running the high-fidelity simulation.
Hence, the offline stage is typically expensive and performed using powerful supercomputers, while the online stage is either real-time or performed with reduced computational effort, both in terms of wall time and resources.

ROMs, as well as the broader class of \emph{surrogate models}, include \emph{data-based}, \emph{equation-based} and \emph{hybrid} models.
In such models, the reduction in computational time is usually obtained by performing a projection of the given parameterized data into a reduced space, where the reduced variables for unseen parameters are either approximated through an interpolation or regression approach (\emph{data-based} approach), or computed by solving a reduced version of the full-order problem (\emph{equation-based} method).

\RALL{The work presented in this manuscript focuses on fully \emph{data-driven}, \emph{non-intrusive} ROMs. Specifically, we focus on models that aim to build a direct mapping from the physical or geometrical parameters of the problem to the output of interest, which are well suited in cases where the governing equations are not available.}
For these models we essentially need a reduction map, its inverse and an interpolation. One of the simplest approaches is to use Proper Orthogonal Decomposition (POD)\cite{lumley1967pod, berkooz1993proper, pearson1901pca, golub1971singular, DeVore2013greedy, ainsworth2021galerkin-nn} to build the reduction map as a linear projection on a smaller space of predetermined (reduced) dimension. POD allows to approximate the solution as a convex combination of the basis spanning the reduced space, namely the \emph{modes}, thus being cheap to compute and interpretable. Nonetheless, it is intrinsically linear, and this may lead to inaccurate approximations especially in the case of scarce data, and highly nonlinear models.

In recent years, many Machine Learning (ML) techniques (such as autoencoders, physics-informed neural networks, generative modelling, graph and convolutional architectures) have been employed to build nonlinear models
\cite{romor2023autoencoder,alla2017nonlinear,amsallem2012nonlinear,kramer2019nonlinear, lee2020model, pichi2024graph, hesthaven2018podnn, coscia2024generative, lee2020model, chen2021physics, hijazi2023pod}, which overcome some limitations of POD while also posing new challenges, one of which being interpretability. Indeed, ML -- especially Deep Learning -- models may be too large and complex to be easily interpreted in the context of the physical problem at hand, and thus provide \emph{black-box} solutions, rather than highlighting physical properties of the system. To keep both the interpretability and the non-linearity, \emph{quadratic models} have been proposed as an extension of POD \cite{xie2018data, geelen2023quadratic-correction, barnett2022quadratic, schwerdtner2024greedy,jain2017quadratic}. In these models, the linear approximation given by POD is enhanced with a quadratic term that still depends on the POD modes. The idea of using known features of the field to model or reconstruct unknown ones is certainly not new in CFD. In fact, it is routinely applied in RANS and LES, where the resolved fields are used to compute \emph{closures} to the governing equations (in terms of Reynolds or sub-grid-scale stresses) \cite{pope2000turbulent,sagaut2006large}. Indeed, since the goal of these quadratic ROMs is to re-introduce the contribution of neglected POD modes, in the form of interactions between the retained ones, we can consider these models a kind of \emph{closure model} for POD. Alternatively, since in the non-intrusive setting the models do not rely on equations, talking about closure models might be 
improper, and thus we may view these models more in general as \emph{corrections} to the base linear ROM, and we can refer to the quadratic term as the \emph{correction term}. 

At the most abstract level, the only prescription about these models is that a quadratic term is added to the POD expansion, in the form of a quadratic operator acting on POD coefficients. However, there is no restriction on how the quadratic operator should be determined. Indeed, a common approach is to infer the best operator from data by solving a linear least-squares problem \cite{geelen2023quadratic-correction,peherstorfer2016quadratic}, which is both computationally efficient and very effective on training data. Models built with this method (which we shall call \emph{Quad-LS}), however, are prone to over-fitting and may not be able to generalise effectively to new unseen data, especially when the least-squares problem is ill-conditioned.

To overcome this issue, we propose a novel approach where the quadratic operator is computed as a nonlinear transformation of the POD modes
. In particular, the key novelties of this work can be summarized as follows:
\begin{itemize}
    \item[(i)] the use of a deep neural network to approximate the quadratic operator. In particular, we employ architectures inspired by the Deep Operator Network (DeepONet) \cite{Lu2021-deeponet} and the Multi-Input Operator Network (MIONet) \cite{jin2022mioinet}, which are the simplest types of neural operators.
    \item[(ii)] use continuous mappings in the spatial coordinates (in the DeepONet approach, named \emph{QuadNet}), where the collocation points are additional inputs of the model;
    \item[(iii)] use continuous mappings in the parameters' space (in the MIONet approach, named \emph{QuadNet-$\mu$}), by adding the parameters as input to the network.
\end{itemize}

As for point (i), the idea is that the increased computational time to learn the best mapping is justified by an improved generalisation capability of the model. Indeed the networks have some intrinsic regularisation that helps to avoid overfitting. We also refer here to another work using DeepONets to learn the residual in non-intrusive ROMs \cite{demo2023deeponet}.

Moreover, we introduce the continuity in space and/or parameters in (ii) and (iii), which is powerful as it makes the model independent of the domain discretisation employed to obtain the high-fidelity training data, and also allows it to be local in parameter space, leading to an improved accuracy. We show that with these methods, that we call \emph{corrected-ROMs}, we get models that 
have a smaller number of parameters than the traditional ones, while also being more accurate. We show that we can train on fewer data, both in terms of spatial sampling points and parameters' instances, 
and still improve the base POD model. This is 
particularly useful for practical applications, where the field data might be available only at few collocation points and/or for a small number of configurations.

The rest of the manuscript is organised as follows. In \cref{sec:methods} we give an overview of non-intrusive ROMs and describe the machine-learning models we employ, in \cref{sec:results} we present the results of numerical experiments on two standard test cases in fluid dynamics, and finally we draw some conclusions in \cref{sec:conclusions}.

\section{Methodology}
\label{sec:methods}

This Section is dedicated to recall the main logic behind non-intrusive ROMs and the machine learning strategies used to enhance and \emph{correct} them.

Non-intrusive ROMs, described in \cref{subsec:met-roms}, are employed to approximate a high-fidelity field of interest for unseen configurations, with significantly reduced computational effort.
The ROMs predictions may be inaccurate in case of data scarcity, especially when employing a linear reduction approach. Recent research works focus on data-driven \emph{closure} strategies specifically designed to mitigate this issue and improve the approximation accuracy \cite{wang2012proper, ahmed2021closures, xie2018data, ivagnes2023pressure, ivagnes2023hybrid}, Such methods are usually applied in intrusive settings and aim to re-integrate the contribution of the discarded basis.
The goal of this work is to extend this framework to non-intrusive settings. In particular, we introduce two different data-driven correction models built using deep operator networks, as extensively described in \cref{subsec:met-corrections}.

\subsection{Non-intrusive ROMs}
\label{subsec:met-roms}
This part is dedicated to briefly recall the theoretical foundations of non-intrusive model order reduction strategies.
Differently from the above-mentioned \emph{intrusive} ROMs, this approach is completely data-driven. Hence, it does not need the knowledge of the mathematical model of the numerical schemes used to obtain the high-fidelity solutions.

Following a classical offline-online pipeline, we collect the high-fidelity solutions, named \emph{snapshots}, corresponding to the field of interest of our problem. \RALL{Let $u \in \mathcal V$ be such field, with $\mathcal V$ an appropriate function space for the problem at hand.} In the numerical results, we focus on the magnitude of the velocity field.
We call $d$ the space dimensionality (we will consider only two-dimensional test cases, hence $d=2$), and we evaluate the velocity field onto the computational mesh points \RALL{$\bm{\hat x} \in \mathbb{R}^{\ndof \times d}$}, being \ndof{} the number of degrees of freedom of the high-fidelity simulations.
 
\RALL{We denote with $\bm{u} = u(\bm{\hat x})$ the vector of nodal values of $u$ at the mesh points $\boldsymbol {\hat x}$ and we consider a parameterized setting with $P$ parameters, namely each snapshot $\bm{u}_i={u}(\bm{\hat x}, \mu_i)$ corresponds to the set of parameters $\mu_i \in \mathbb{R}^P$.}
After collecting the snapshots, we can build the snapshots' matrix:
\RALL{\begin{equation*}
    \mathbf{S}=\begin{bmatrix}
\vert&\vert&&\vert\\
\bm{u}_1 &\bm{u}_2 &\dots&\bm{u}_{\ntrain}\\
\vert&\vert&&\vert
    \end{bmatrix}.
\end{equation*}}
Such matrix has dimension $\ndof \times \ntrain$, where \ntrain{} is the number of snapshots considered in the offline stage.

Non-intrusive ROMs consist of three main steps:
\begin{itemize}
    \item \textbf{A reduction step }$\mathcal{R}$: projection of the solutions' manifold into a space of reduced dimensionality $r \ll \ndof$.
    \item \textbf{An approximation step }$\mathcal{A}$: approximation of the solution for unseen parameters in the reduced space through interpolation or regression methods.
    \item \textbf{Backmapping}: the solution is backmapped into the space with full dimensionality.
\end{itemize}

The reduction step $\mathcal{R}$ is here achieved with a Proper Orthogonal Decomposition (POD) technique, which performs a linear projection of the snapshots through a Singular Value Decomposition (SVD) of matrix $\mathbf{S}$.
The main POD hypothesis is that the velocity field can be approximated as a convex combination of the basis spanning the reduced space, namely the \emph{modes}.
Computing the SVD of matrix $\mathbf{S}$, we indeed obtain three matrices \RALL{$\bm \Phi$, $\bm \Sigma$, and $\bm V$} such that:
\begin{equation*}
\RALL{\mathbf{S}=\bm \Phi \bm \Sigma \bm V^T}.
\end{equation*}

Typically, a reduced dimension $r$ is \emph{a-priori} selected depending on the corresponding retained energy. 
\RALL{The first $r$ columns of $\bm \Phi$ are the vectors $\bm \phi_i \in \mathbb R^{\ndof},\; i= 1,\dots,r$, corresponding to the nodal values of the modes $\phi_i$ in the mesh points $\boldsymbol {\hat x}$, i.e. $\bm \phi_i = \phi_i(\bm{\hat x})$.}
The velocity can then be approximated as follows:
\begin{equation}
    \RALL{{u}(\bm{\hat x}, \mu) \simeq \tilde{{u}}(\bm{\hat x}, \mu)= \sum_{i=1}^r  a_i(\mu) \phi_i(\bm{\hat x}),}
    \label{eq:rom-approx}
    \end{equation}
where $\{a_i\}_{i=1}^r$ are the \emph{reduced coefficients} or \emph{variables}.
The computational gain of the POD is highlighted in the expression \cref{eq:rom-approx} by the separation of the \textbf{space} dependency (embedded in the modes) and the \textbf{parametric} dependency (in the reduced variables).

We specify that the coefficients for the train parameters are directly computed as:
\begin{equation*}
    \RALL{\bm{a}(\mu_j)=\Phi_r^T \bm{u}_j(\mu_j)}
 \quad j=1, \dots, \ntrain.
 \end{equation*}

Since we may want to employ the ROM to infer the velocity in unseen configurations $\mu^{\star}$, the key problem is now how to build a map $\mathcal{A}$ such that:
$$\bm{a}(\mu^{\star})=\mathcal{A}(\mu^{\star})=\{a_i(\mu^{}\star)\}_{i=1}^r.$$

Different mappings can be employed, such as the Radial Basis Function (RBF) interpolation~\cite{acosta1995radial, buhmann2000radial}, the Gaussian Process Regression (GPR)~\cite{chu2005gaussian, ortali2022gaussian}, or other types of regression like deep neural networks~\cite{hesthaven2018podnn}.
In this work, we employ a RBF interpolation, where:
\begin{equation}
\mathcal{A}(\mu^{\star})=\sum_{i=1}^{\ntrain}\bm{\omega}_i \varphi(\|\mu^{\star} - \mu_i\|),
    \label{eq:rbf}
\end{equation}
where $\varphi(\cdot)$ is the radial basis functions' kernel, having weights $\bm{\omega}_i$, and centers $\mu_i$, $i=1, \dots, \ntrain$.
The values of the weights are simply computed from the conditions  
\begin{equation*}
\bm{a}(\mu_j)=\sum_{i=1}^{\ntrain}\bm{\omega}_i \varphi(\|\mu_j - \mu_i\|), \quad j=1, \dots, \ntrain.
    \label{eq:rbf-data}
\end{equation*}

The RBF interpolation allows to obtain efficient and precise approximations, and to work with unstructured data. Moreover, it allows flexibility in the kernel $\varphi(\cdot)$ choice. It is usually selected depending on the specific test case and on the complexity of the reduced manifold.
In our cases, we choose a \emph{thin plate spline} and a \emph{linear} kernel, namely:
\begin{equation}
\varphi(d)=d^2 \log{(d)} ~~\text{ and }~~ \varphi (d) = d.
    \label{eq:rbf-kernel}
\end{equation}

Once the reduction and approximation steps are performed, the reduced variables $\bm{a}(\mu^{\star})$ can be backmapped using expression \cref{eq:rom-approx}.

\subsection{Corrected-ROMs via deep operator networks}
\label{subsec:met-corrections}

As already specified, the reconstruction in \cref{eq:rom-approx} may not be accurate when keeping a small dimension $r$, especially in data scarcity regimes. This motivates the introduction of \RB{correction} terms, namely:
\begin{equation}
\RALL{{u}(\bm{\hat x}, \mu)\simeq \sum_{i=1}^r a_i(\mu) \phi_i(\bm{\hat x}) \, + \, {\tau}(\bm{\hat x}; \mu) .}
    \label{eq:rom-closed}
\end{equation}

The correction term \RALL{$\bm\tau (\mu)= {\tau}(\bm{\hat x}; \mu)$} is modelled through a mapping that we call $\mathcal{M}$, and it is completely \emph{data-driven}. The goal is indeed to approximate the \emph{exact} correction:
\begin{equation}
    \RALL{\bm \tau ^{exact} (\mu) = {\tau}^{exact}(\bm{\hat x};\mu)={u}(\bm{\hat x}, \mu)-\tilde{{u}}(\bm{\hat x}, \mu).}
    \label{eq:corr-exact}
\end{equation}

The equation \cref{eq:corr-exact} aims intuitively at \emph{correcting} the POD approximation reintegrating the data knowledge.
The exact correction is modeled with a mapping, that should in general depend on the known quantities, namely the modes, the reduced variables, and, eventually, the parameters: $\bm{\tau}(\phi_1, \dots, \phi_r, a_1,\dots, a_r; \mu)$.
In particular, we compare the performance of the following maps, all providing a quadratic ansatz in the reduced coefficients:

\begin{itemize}
    \item \textbf{Quad-LS} a non-trainable approach, \RB{already existing in literature \cite{geelen2023quadratic-correction}} and briefly recalled in \cref{subsec:met-quad-ls};
    \RALL{\item \textbf{QuadNet}: a DeepONet-inspired architecture having inputs $(\bm \phi_1, \dots, \bm \phi_r)$ and $ (a_1,\dots, a_r)$, described in \cref{subsec:met-quadnet};
    \item \textbf{QuadNet-$\mu$}: a MIONet-inspired network having as inputs $(\bm \phi_1, \dots, \bm \phi_r)$, $ (a_1,\dots, a_r)$ and $\mu$. The method is described in \cref{subsec:met-quadnet-mu}. }
\end{itemize}

\subsubsection{Quadratic Least-Squares model: Quad-LS (state-of-the-art)}
\label{subsec:met-quad-ls}
\hfill\\
The Quad-LS strategy aims at approximating the correction term with the quadratic form:
\begin{equation}
    \RALL{\bm{\tau}(\mu)=\bm{a}^T(\mu) \, \boldsymbol C \, \bm{a}(\mu).}
    \label{eq:map-quad-ls}
\end{equation}
\RB{This strategy is employed in \cite{geelen2023quadratic-correction}, and relies on the determination of the quadratic operator $\bm C$ via least-squares optimisation. In fact, the procedure to determine $\bm C$ is very general, and can be applied to the inference of quadratic operators in other contexts as well, as is done for instance in \cite{peherstorfer2016quadratic}}.

\RALL{We stress that in \cref{eq:map-quad-ls}, $\bm \tau(\mu) = \tau (\bm{\hat x};\mu) \in \mathbb R^{\ndof}$ is the $\ndof$--dimensional vector with the values of the correction at each mesh node, while $\boldsymbol C\in \mathbb{R}^{r \times \ndof \times r}$ is a tensor that satisfies the following minimisation problem:
\begin{equation}
    \min_{{\boldsymbol C}\in\mathbb{R}^{r \times \ndof \times r}} \sum_{j=1}^{\ntrain}\|{\tau}^{exact}(\mu_j)- \bm{a}(\mu_j)^T \, \boldsymbol C \, \bm{a}(\mu_j)\|_2^2.
    \label{eq:min-ls}
\end{equation}}

The minimization problem in \cref{eq:min-ls} can be re-written and solved as a least squares problem exploiting the symmetry of the quadratic form $\bm{a}^T \,\RALL{\boldsymbol C} \, \bm{a}$.
We can indeed write the quadratic form as a double summation:
\begin{equation}
    \sum_{i=1}^r \sum_{j=1}^r a_i(\mu)a_j(\mu)C_{ij}.
\end{equation}
We introduce a vector \RA{$\tilde{\bm{a}}(\mu) \in \mathbb{R}^\mathcal S$} defined as 

\begin{equation}
    \tilde{\bm{a}}(\mu) = \begin{bmatrix}
        \bm{a}^{(1)}(\mu)\\
        \vdots \\
        \bm{a}^{(r)}(\mu)
    \end{bmatrix}, \text{where }
    \bm{a}^{(i)}(\mu) = a_i(\mu) \begin{bmatrix}
        a_1(\mu)\\
        \vdots \\
        a_i(\mu)
    \end{bmatrix} \quad \in \mathbb{R}^i.
\end{equation}
The vector $\tilde{\bm{a}}(\mu)$ contains all pairwise products of the components of $\bm{a}$ with themselves, without repetition. Hence it has dimension $\mathcal{S}=\frac{r(r+1)}{2}$.

We can also write \RALL{$\boldsymbol C $} in block form, as
\begin{equation}
    \RALL{\hat{\boldsymbol C} = \begin{bmatrix}
        \hat{\boldsymbol C}^{(1)}, \dots, \hat{\boldsymbol C}^{(r)}
    \end{bmatrix}} \quad \in \mathbb{R}^{\ndof \times \mathcal{S}},
    \label{eq:block-C}
\end{equation}
where $\RALL{\hat{\boldsymbol C}^{(i)}} \in \mathbb{R}^{\ndof \times i}$. We can finally rewrite the minimization in \cref{eq:min-ls} as:
\begin{equation}
    \RALL{\min_{\hat{\boldsymbol C}\in\mathbb{R}^{\ndof\times \mathcal{S}}} \sum_{j=1}^{\ntrain} ||\boldsymbol{\tau}(\mu_j) - \hat{\boldsymbol C}\tilde{\bm{a}}(\mu_j)||_2^2,}
\end{equation}

We refer the reader to \cite{peherstorfer2016quadratic} for the extended least squares derivation.

It is worth remarking that this approach allows to find a \emph{unique} tensor \RALL{$\boldsymbol C$}. As we will see in the numerical results, the fact that \RALL{$\boldsymbol C$} is not parameter-dependent may lead to inaccurate test approximations, especially when the training manifold (used to solve \cref{eq:min-ls}) does not reflect the real complexity of the high-fidelity model.

\subsubsection{Quadratic Network model: QuadNet}
\label{subsec:met-quadnet}
\begin{figure}[!htb]
    \centering
    \includegraphics[width=0.8\linewidth]{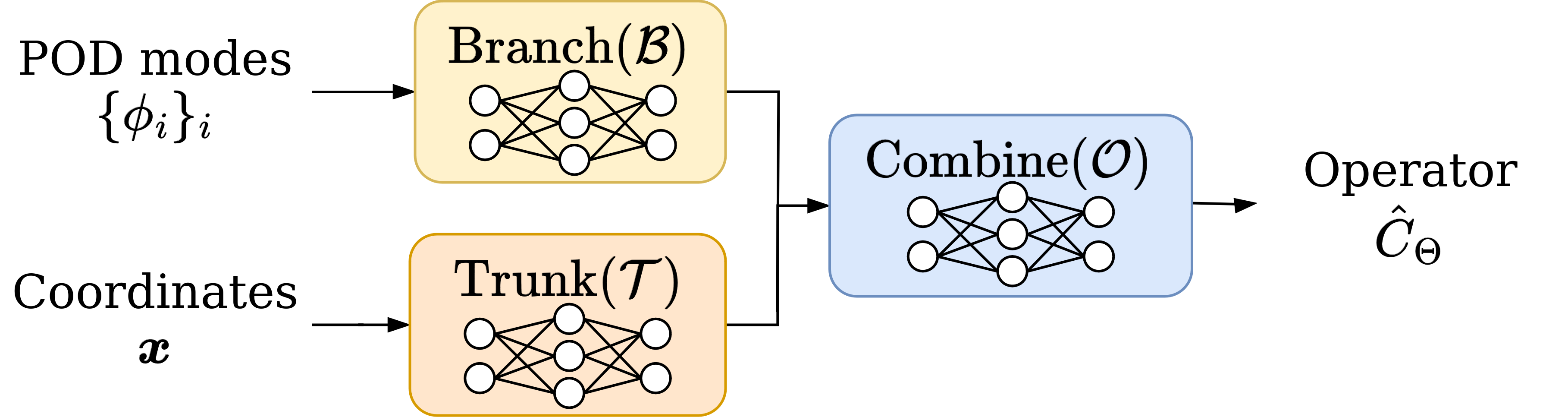}
    \caption{Schematic representation of the DeepONet architecture employed in \quadnet. }
    \label{fig:deeponet-sketch}
\end{figure}

To improve the generalization capability of the model described in the previous \cref{subsec:met-quad-ls}, we provide in this work a novel machine-learning based alternative. In particular, we propose the following correction model:
\begin{equation}
     \bm{\tau}(\bm{x}; \mu)=\bm{a}^T(\mu)C(\Phi_r, \bm{x})\bm{a}(\mu),
    \label{eq:map-quadnet}
\end{equation}
\RALL{where, the dependence on $\bm x$ highlights that the correction is now a \emph{continuous} function of space, and tensor $C(\Phi_r, \bm{x})$ is modelled through a DeepONet-inspired architecture, displayed in \cref{fig:deeponet-sketch}}. The DeepONet is one of the simplest examples of neural operator for its peculiar architecture. It is indeed composed of two or more \emph{sub-networks}, that separately handle different inputs.

In our case \RALL{-- since we are interested in the corrections at the mesh nodes --} the modes matrix \RB{$\Phi_r=[\bm \phi_1, \dots, \bm \phi_r] \in \mathbb{R}^{\ndof \times r}$ is the input to the \emph{branch network} $\mathcal{B}$, while the mesh points $\bm{\hat x} \in \mathbb{R}^{\ndof \times d}$ are the inputs to the \emph{trunk network} $\mathcal{T}$}. More precisely, we define $\mathcal B : \mathbb{R}^r \mapsto \mathbb{R}^p $ and $\mathcal T : \mathbb{R}^d \mapsto \mathbb{R}^p $, so that a single input to the DeepONet consists of an $r$-dimensional vector, representing the $r$ modes evaluated at a given point in the domain, and a $d$-dimensional vector, representing the coordinates of said point. 
The outputs of the two networks are then combined together and processed by a third sub-network $\mathcal{O}:\mathbb R^p \mapsto \mathbb{R}^\mathcal{S}$, whose output is the final matrix $\hat{C}_{\Theta}(\Phi_r, \boldsymbol{x})$, namely the blocks expression in \cref{eq:block-C}, where $\Theta$ are the network parameters. Note that although the dimension $p$ can be set to an arbitrarily large value, from experiments we found that a suitable value is $p = \mathcal S$.


The novelty of this model is that it is \emph{space continuous}. On the one hand, the quadratic operator can be evaluated in points not belonging to the original mesh. On the other hand, the DeepONet may be trained on a subset of the original mesh, limiting the computational burden and the training times.
The numerical results also include a sensitivity analysis on such training size.

The loss function we try to minimise here can be written as
\begin{equation}
\label{eq:quadnet-loss}
    \mathcal{L}(\Theta) = \frac{1}{N_\mu}\sum_{i=1}^{N_\mu} \frac{ \lVert \boldsymbol{\tau}_i - (\boldsymbol{a}^{(2)}_i)^T \hat{C}_{\Theta}(\Phi_r, \RALL{\bm{\hat x}})\rVert^2_{l^2} }{\lVert \boldsymbol{\tau}_i\rVert^2_{l^2}} ,
\end{equation}
where $\boldsymbol{\tau}_i$ is the exact correction evaluated on the mesh, $(\boldsymbol{a}^{(2)}_i)^T$ is the $\mathcal{S}$-dimensional vector of pairwise products between modal coefficients, $\hat{C}_{\Theta}$ is the quadratic operator evaluated on the mesh (output of the DeepONet).

\remark{
The choice of a \emph{relative error} loss over the more common mean squared error loss is motivated by the necessity of making the training less sensitive to different scales in the correction terms. Indeed, from experiments it was clear that the training procedure is quite delicate and sometimes sensitive to the choice of hyperparameters such as the learning rate. This happens especially when the corrections have 
small values ($\mathcal{O}(10^{-3})$) or 
have different orders of magnitude. We remark that, even though the normalisation terms can be relatively expensive to compute for large meshes, they can all be computed prior to the training and thus do not represent a significant computational burden.}

\subsubsection{\texorpdfstring{Quadratic Parametric Network model: QuadNet-${\mu}$}{Quadratic Parametric Network model}}
\label{subsec:met-quadnet-mu}
The model presented in the previous parts \cref{subsec:met-quad-ls} and \cref{subsec:met-quadnet} are fixed in the parametric space.
We propose in this part a further parameterized extension of the model presented in \cref{subsec:met-quadnet}. The novel QuadNet-$\mu$ model approximates the correction as:
\begin{equation}
     \RALL{{\tau}(\bm{x}; \mu)=\bm{a}^T(\mu)C(\Phi_r, \bm{x}, \mu)\bm{a}(\mu).}
    \label{eq:map-quadnet-mu}
\end{equation}

In this case, the neural network used to compute $C(\Phi_r, \bm{x}, \mu)$ has a MIONet-inspired structure, which is the extension of the DeepONet with more than two sub-networks.
Indeed, the MIONet is composed of two branch networks and one trunk network, that separately process the POD modes $\Phi_r$, the parameter $\mu$, and the spatial coordinates $\boldsymbol{x}$. More precisely, the first branch network defines a map $\mathcal{B}_1 : \mathbb R ^r\mapsto \mathbb R^p$, whose input is an $r$-dimensional vector representing the $r$ modes evaluated at a given point in the domain; the second branch network is a map $\mathcal{B}_2: \mathbb R^{d_\mu} \mapsto \mathbb R^p$, whose input is a $d_\mu$-dimensional vector representing the parameters of the problem (in our case $d_\mu=1$); finally, the trunk network is a map $\mathcal{T}: \mathbb R ^d\mapsto \mathbb R^p$, whose input is the $d$-dimensional vector of spatial coordinates for a given point in the domain. As for the \quadnet, we combine the outputs of these networks with a fourth sub-network $\mathcal{O}:\mathbb R^p \mapsto \mathbb{R}^\mathcal{S}$, which outputs the block version of operator $C$, namely $\hat{C}_{\Theta}(\Phi_r, \boldsymbol{x}, \mu)$, where $\Theta$ are the parameters of the network learned during training.

\begin{figure}[!htb]
    \centering
    \includegraphics[width=0.8\linewidth]{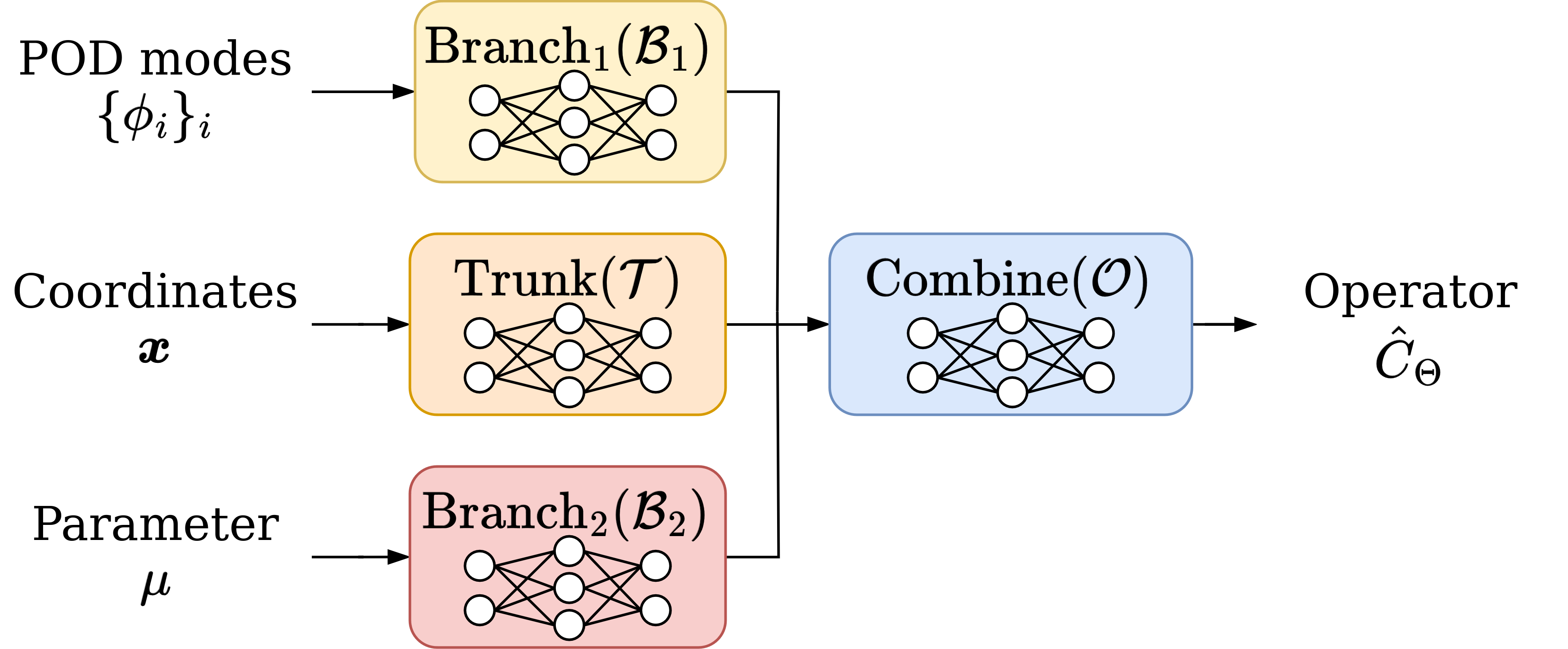}
    \caption{Schematic representation of the MIONet architecture employed in \quadnetmu.}
    \label{fig:mionet-sketch}
\end{figure}

 The architecture employed in this case is represented in \cref{fig:mionet-sketch}.


The loss function minimized during training is the same as in \cref{eq:quadnet-loss}, but $\hat{C}_{\Theta}=\hat{C}_{\Theta}(\Phi_r, \boldsymbol{x}, \mu)$.
\section{Results}
\label{sec:results}
We assess the quality of our models by applying them to two standard benchmarks in fluid dynamics, the backward-facing step and the lid-driven cavity flows. 
This part is structured as follows:
\begin{itemize}
\item We shall begin by describing the test problems and their full-order discretisations, in \cref{subsec:fom-1} and \cref{subsec:fom-2}, for the two test cases, respectively.
\item A description of the hyperparameters of the different models is provided in \cref{subsec:res-roms};
\item We proceed with the comparison among the novel methods and the linear and quadratic baseline models (\cref{subsec:comparison}), considering fixed number of modes $r$ and of snapshots $N_{\mu}$;
\item We also investigate the accuracy of our model in two different data regimes:
\begin{itemize}
    \item \emph{partial} data (\cref{subsec:partial}), training the models with only a subset of the mesh nodes, keeping the values of $r$ and $N_{\mu}$;
    \item \emph{scarce} data (\cref{subsec:scarce}), training the networks considering different parameters' instances $N_{\mu}$, for different values of $r$.
\end{itemize}
\end{itemize}

\subsection{Case 1: backward-facing step}
\label{subsec:fom-1}
\begin{figure}
    \centering
    \begin{subfigure}{0.55\textwidth}
    \centering
        \includegraphics[width=\textwidth]{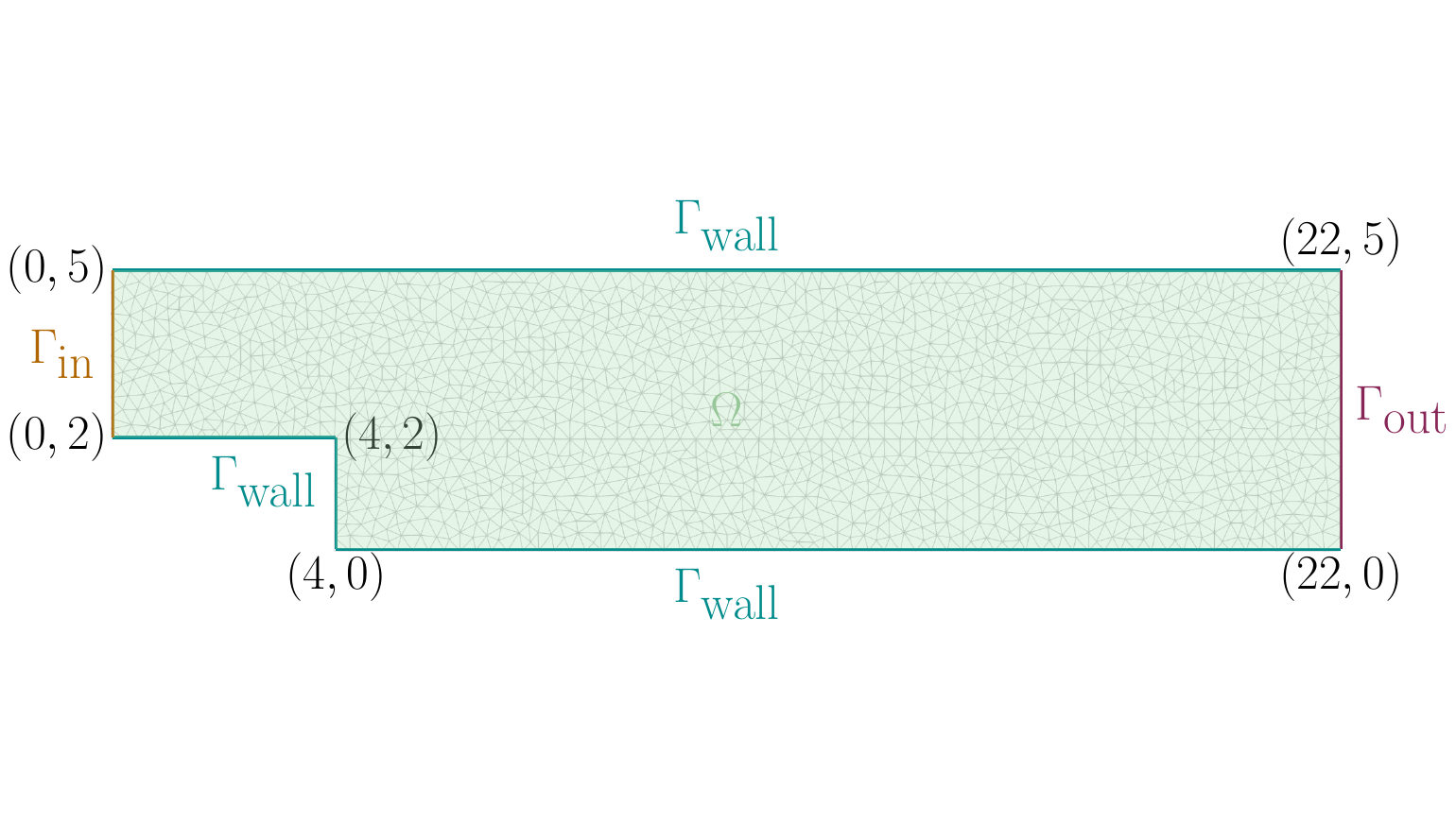}
        \caption{Backward-facing step.}
        \label{fig:backstep-domain}
    \end{subfigure}
    \begin{subfigure}{0.35\textwidth}
    \centering
        \includegraphics[width=\textwidth]{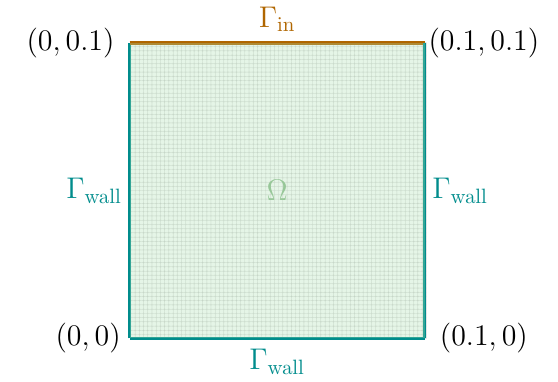}
        \caption{Lid-driven cavity.}
        \label{fig:cavity-domain}
    \end{subfigure}
    \caption{Domains and corresponding offline grid for the two test cases.}
\end{figure}

In this case we wish to find steady-state solutions to the incompressible Navier-Stokes equations in the domain represented in \cref{fig:backstep-domain}. The equations are complemented with no-slip conditions at the top and bottom boundaries, homogeneous Neumann conditions for the Cauchy stress tensor at the outlet, and Dirichlet conditions at the inlet, which impose a parabolic Poiseuille profile for the velocity. The full system of equations is given by
\begin{equation} 
\label{eq:backstep-governing}
\begin{cases}
    \frac{\partial u_i}{\partial x_i} = 0 \quad \text{in}~\Omega,\\  
    u_j \frac{\partial u_i}{\partial x_j} = -\frac{\partial p}{\partial x_i} + \nu \frac{\partial^2 u_i}{\partial x_j \partial x_j} \quad \text{in}~\Omega ,\\
     \RALL{u_0 = u_1 = 0} \quad \text{on} ~\Gamma_{wall}, \\ 
    \RA{ n_j \left ( -p\delta_{ij} + \nu  \frac{\partial u_i}{\partial x_j} \right ) = 0  \quad \text{on} ~\Gamma_{out},} \\
     u_0 = \frac{\mu}{2.25} (x_1 -2)(5-x_1) \quad \text{on} ~\Gamma_{in}, \\
    u_1 = 0  \quad \text{on} ~\Gamma_{in} ,
\end{cases}
\end{equation}
\RA{where $p$ and $\nu$ are the kinematic pressure ($p = P/\rho$, with $P$ dynamic pressure) and viscosity, $\delta_{ij}$ is the Kronecker symbol, $\boldsymbol{n} = (n_0,n_1) $ is the outward normal unit vector to $\Gamma_{out}$, and $\mu \in [1,80]$ is a parameter which gives the amplitude of the parabolic profile}.

The full-order solution is computed via the Finite Element method, implemented in the open-source Python package FEniCS \cite{baratta2023dolfinx}. Problem \cref{eq:backstep-governing} is written in weak form as:

\emph{Given $\mu$, find $\mathbf{u} \in \mathcal{V}$, $p\in \mathcal{Q}$ such that the following equations are satisfied:}
\RA{
\begin{equation}\label{eq:backstep-weak}
    \begin{cases}
        \nu \int_{\Omega} \nabla \mathbf{u} : \nabla \mathbf{v} ~d\Omega + \int_{\Omega} [(\mathbf{u} \cdot \nabla)\mathbf{u}] \cdot \mathbf{v} ~d\Omega - \int_{\Omega} p ~ \nabla\cdot\mathbf{v} ~d\Omega  = 0  \quad \forall ~\mathbf{v}\in\mathcal{V}, \\
        \int_{\Omega} q~\nabla\cdot\mathbf{u} ~d\Omega  = 0 \quad \forall ~q\in\mathcal{Q},
    \end{cases}
\end{equation}
}
which can be discretised by approximating the 
infinite-dimensional spaces $\mathcal{V}, \mathcal{Q}$ with finite-dimensional Taylor-Hood $\mathbb{P}^2-\mathbb{P}^1$ polynomial spaces. 
The domain is discretised with 3091 triangular elements, and the discrete problem is solved for $N_{\mu}=500$ values of $\mu$. For each $\mu$ we save the values of the velocity and pressure fields at the 1639 vertices of the triangulation.

\subsection{Case 2: lid-driven cavity}
\label{subsec:fom-2}
In this case we aim at predicting the fully developed flow fields for incompressible fluids in the domain represented in \cref{fig:cavity-domain}. As for the previous test case, the governing equations are the Navier-Stokes equations, now in their time-dependent form. The boundary conditions here are non-slip conditions on $\Gamma_{wall}$ and Dirichlet conditions on $\Gamma_{in}$, where we impose a constant, horizontal velocity, whose \RA{magnitude} is given by a parameter. The equations for this test case read as follows:
\begin{equation} 
\label{eq:cavity-governing}
\begin{cases}
    \frac{\partial u_i}{\partial x_i} = 0 \quad \text{in}~\Omega\times\mathcal{I},\\  
    \frac{\partial u_i}{\partial t} + u_j \frac{\partial u_i}{\partial x_j} = -\frac{\partial p}{\partial x_i} + \nu \frac{\partial^2 u_i}{\partial x_j \partial x_j} \quad \text{in}~\Omega\times \mathcal{I} ,\\
    \RALL{u_0 = u_1 = 0} \quad \text{in} ~ \Omega\times \{0\}, \\
    \RALL{u_0 = u_1 = 0} \quad \text{on} ~ \Gamma_{wall}\times\mathcal{I}, \\
    u_0 = \mu \quad \text{on} ~ \Gamma_{in}\times\mathcal{I}, \\
    u_1 = 0 \quad \text{on} ~ \Gamma_{in}\times\mathcal{I} ,
\end{cases}
\end{equation}
where $\mathcal{I} = [0,5]\, \si{\second}$ is the time domain, $\nu = 10^{-5}\, \si{\meter\squared\per\second}$ \RA{and $p$ are the kinematic viscosity and pressure}, and $\mu \in [0.5,1]\, \si{\meter\per\second}$ is the velocity magnitude at the top boundary. The Reynolds number for this problem can be defined as \RA{$Re=\frac{\mu L}{\nu}$, where $L=0.1\, \si{\meter}$} is the width (and height) of the domain. We stress that, for this choice of parameters, the Reynolds number varies in the range $[5\times10^3,10^4]$, \RA{meaning that the flow is in a transitional state from laminar to turbulent, and it is dominated by convection}.

The full-order model consists in a Finite Volume discretisation of the 
Navier-Stokes equations, implemented in the open-source software OpenFOAM \cite{jasak1996error,moukalled2016finite,schreiber1983driven} . The spatial domain is discretised with a uniform grid of $70\times70\times1$ hexahedral cells, and the time domain is discretised at steps $\delta t=\si{\num{1e-4}\, \second}$. A dataset was built by solving the problem for 150 values of $\mu$ and saving a snapshot of the fields evaluated at 5041 points in the domain after $5\si{\second}$.

\subsection{Setting of Corrected-ROMs}
\label{subsec:res-roms}

Here we provide details of the models described in \cref{sec:methods}, as implemented for the two test cases. The models are trained on $N_{\mu} = 400$ and $N_{\mu} = 100$ high-fidelity snapshots for the backward-facing step and lid-driven cavity flows, respectively. The remaining snapshots are employed for testing only.

\subsubsection{POD-RBF}
This is the linear model described in \cref{subsec:met-roms}, which employs POD in the reduction step and RBF interpolation. The main parameter we need to fix for this model is the number $r$ of POD modes we want to retain. Since we are mainly interested in exploring the under-resolved regime, we chose $r=3$ for both test cases. A second parameter that can determine the performance of this model is the kernel function employed for the RBF interpolation. Among the different possibilities, we settled for a linear kernel for the backward-facing step and a thin plate spline kernel for the lid-driven cavity (see \cref{eq:rbf-kernel}).

\begin{figure}[!htb]
    \centering
    \includegraphics[width=0.8\linewidth]{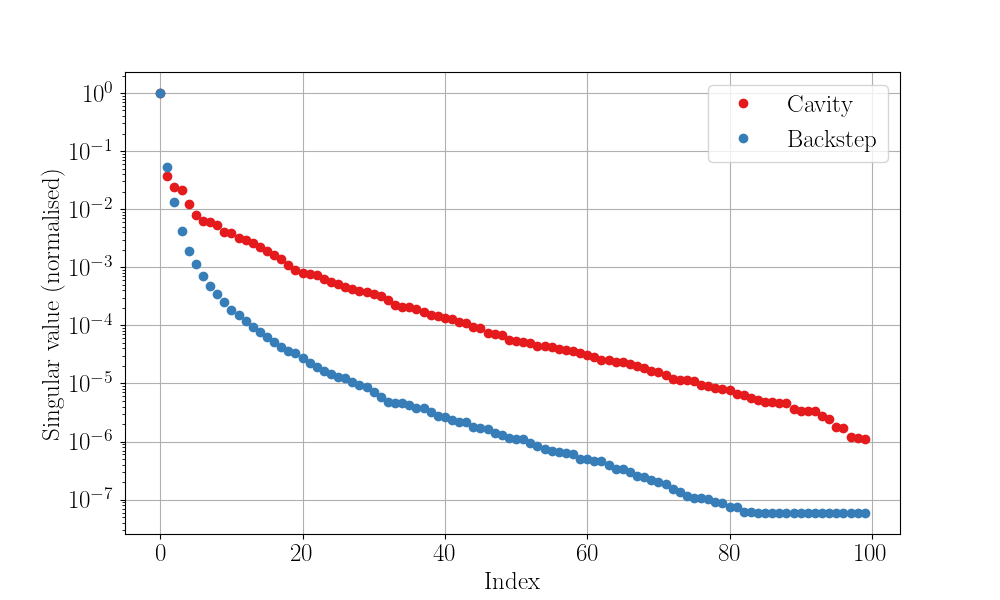}
    \caption{Decay of the snapshot matrix's singular values for the two test cases. Values are normalised to the first (largest) singular value.}
    \label{fig:singular-values}
\end{figure}

In \cref{fig:singular-values} we show the normalised singular values of the snapshot matrix for the two test cases. The decay of the singular values suggests that \podrbf may be better suited for the backward-facing step than for the lid-driven cavity. This is something we expected a priori, since we know that the cavity flow is \RA{convection dominated} and hence highly nonlinear.

\subsubsection{Quad-LS}
For this model there are no major parameters to fix, as choosing $r$ automatically determines the size $\mathcal{S}$ of the quadratic operator. In particular, for our choice $r=3$ we get $\mathcal{S}=6$. We remind the reader that with this method we find a tensor that represents the $\mathcal{S}$ entries of the operator evaluated on the $\ndof$ discretisation points in the domain. Although here we have $\ndof = 1639 $ for the backward-facing step and $\ndof = 5041 $ for the lid-driven cavity, for more complex problems and more POD modes the size of the tensor $\hat C$ can grow very rapidly. In these cases, the minimisation problem \cref{eq:min-ls} can be easily solved with the least squares solver implemented in PyTorch.

\subsubsection{QuadNet}
For this model we need to decide on a suitable architecture for the DeepONet. The simplest choice is to employ fully connected feed-forward networks both for branch and trunk. As described in \cref{subsec:met-quadnet}, we also employ a third network to combine the output of branch and trunk in a nonlinear fashion. This small network is also a feed-forward. 

A summary of the hyperparameters' setting for this model can be found in \cref{tab:deeponet-arch}. For both test cases these architectures seem effective and the only minor differences may be on the training hyperparameters.
\begin{table}[!htb]
\centering
    \caption{Summary of the network parameters for the \quadnet model.}
    \begin{tabular}{cccccccc}\toprule
    &\multicolumn{4}{c}{Architecture}&\multicolumn{3}{c}{Training}\\
    \cmidrule(lr){2-5} \cmidrule(lr){6-8}
    & {\small Input} & {\small Hidden} & {\small Output} & 
    {\small Activation} & {\small lr} & {\small Min. loss} & {\small Max. epochs}\\\midrule
    $\mathcal{B}$ & 3 & $7\times20$ & 6 & Tanh & & & \\
    $\mathcal{T}$ & 2 & $7\times20$ & 6 & Tanh & $O(10^{-3})$ & $10^{-2}$ & 20000 \\
    $\mathcal{O}$ & 6 & / & 6 & Tanh & & & \\
    \bottomrule
    \end{tabular}
    \label{tab:deeponet-arch}
\end{table}

\subsubsection{\texorpdfstring{QuadNet-$\mu$}{QuadNet (parametric)}}
This model is essentially the same as \quadnet, only with the addition of a second branch network that processes the parameter $\mu$. A summary of model hyperparameters is given in \cref{tab:mionet-arch}.
\begin{table}[!htb]
\centering
    \caption{Summary of the network parameters for the \quadnetmu model.}
    \begin{tabular}{cccccccc}\toprule
    &\multicolumn{4}{c}{Architecture}&\multicolumn{3}{c}{Training}\\
    \cmidrule(lr){2-5} \cmidrule(lr){6-8}
    & {\small Input} & {\small Hidden} & {\small Output} & 
    {\small Activation} & {\small lr} & {\small Min. loss} & {\small Max. epochs}\\\midrule
    $\mathcal{B}_1$ & 3 & $7\times20$ & 6 & Tanh & \multirow{4}{*}{$O(10^{-3})$} & \multirow{4}{*}{$10^{-2}$} & \multirow{4}{*}{20000}\\
    $\mathcal{B}_2$ & 1 & $7\times20$ & 6 & Tanh & & & \\
    $\mathcal{T}$ & 2 & $7\times20$ & 6 & Tanh & & & \\
    $\mathcal{O}$ & 6 & / & 6 & Tanh & & & \\
    \bottomrule
    \end{tabular}
    \label{tab:mionet-arch}
\end{table}

This model is also trained by minimising the loss function \cref{eq:quadnet-loss}. The only difference in training with respect to \quadnet is that this model is slightly larger and more complex, and as such it could benefit from smaller learning rates.

At this point we can make a first comparison between the three quadratic models. We already noted that one of the major differences between \quadls and the NN-based models is that the former is discrete in space, while the latter are continuous. We have seen that this implies that in \quadls the quadratic operator is determined by $\mathcal S \times\ndof$ parameters, which in these cases amounts to $\approx10000$ for the backward-facing step and $\approx 30000$ for the lid-driven cavity. At the same time, from \cref{tab:deeponet-arch,tab:mionet-arch} we can compute that \quadnet and \quadnetmu have $\approx 5700$ and $\approx 8400$ parameters respectively, making them effectively smaller than \quadls.
\RALL{On the one hand, having a smaller model is advantageous in terms of memory utilisation, especially for large-scale real-world problems where meshes are substantially larger than the ones considered in this work. On the other hand, for these test cases, we could not observe any benefit in terms of computational time related to the reduced size of NN-based models. In fact, while the online cost is comparable for the three models, the offline phase for NN-based ones is much more expensive compared to that of \quadls. The typical training times for \quadnet and \quadnetmu -- on the machine considered for these experiments\footnote{\RALL{The machine is a workstation equipped with an 11th Gen Intel(R) Core(TM) i7-11700 CPU, 32GB of RAM and an NVIDIA Quadro RTX 4000 GPU with 8GB of memory.}} -- range from $\approx 8$ to $\approx 25$ minutes (depending on the number of parameters and the training set size), while \quadls requires just a few milliseconds (typically less than $10$) to train. This increased offline cost is readily justified when put in perspective with the time required to compute the full-order solution. Consider for example the lid-driven cavity flow described earlier. For that model, a single full-order run takes about $3$ minutes to complete, which means that in the time required to train the network we could get roughly $3$ to $8$ more FOM solutions, depending on the neural network size. This means that if we train both \quadls and \quadnetmu on the same number of snapshots, in order to amortise the additional cost required by \quadnetmu we would just need to query the model $3$ to $8$ more times with respect to \quadls, which is likely a small amount in a real-world many-query context. }

\subsection{Comparison}
\label{subsec:comparison}
\begin{figure}[!htb]
    \centering
    \begin{subfigure}{\linewidth}
        \centering
        \includegraphics[width=\linewidth, trim={5cm 0 6cm 0},clip]{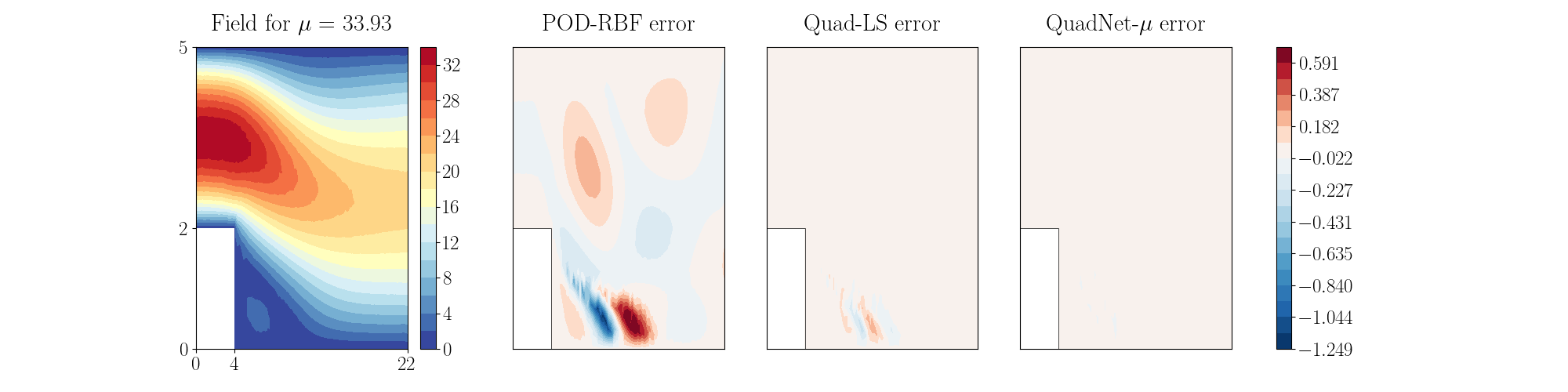}
        \caption{Backward-facing step.}
        \label{fig:backstep-errors}
    \end{subfigure}
    \begin{subfigure}{\linewidth}
        \centering
        \includegraphics[width=\linewidth, trim={5cm 0 6cm 0},clip]{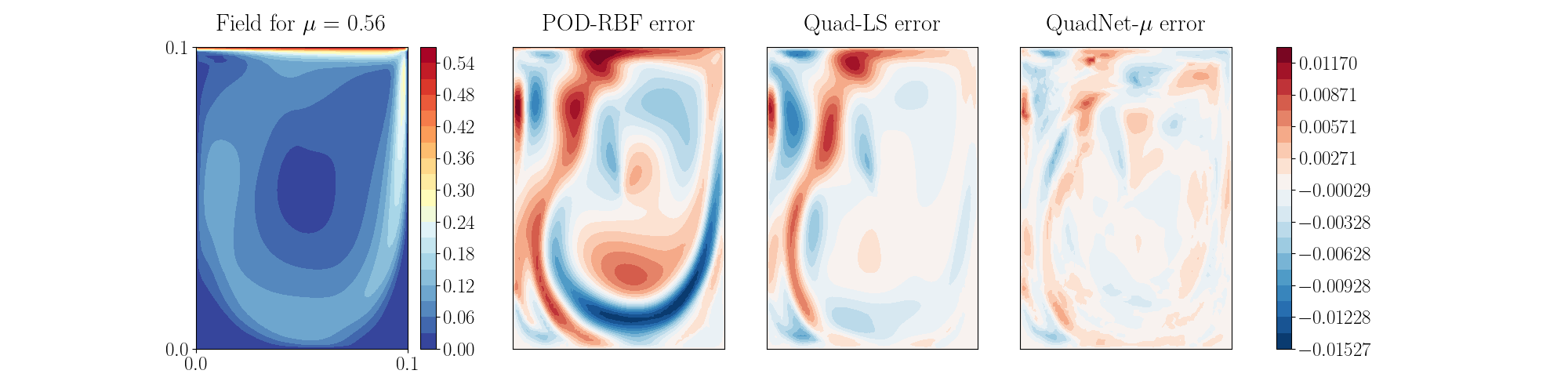}
        \caption{Lid-driven cavity.}
        \label{fig:cavity-errors}
    \end{subfigure}
    \caption{Plots for a random test FOM solution and corresponding error fields for \podrbf, \quadls and \quadnetmu. Here $\mu = 33.93$ for the backward-facing step and $\mu = 0.56$ for the lid-driven cavity.}
    \label{fig:error-fields}
\end{figure}

We start by comparing the four models in terms of their accuracy. We train all four models on the same data for the two test cases, and compute the relative error on the prediction for both train and test datasets. Specifically, given the FOM solution $\boldsymbol{u}(\mu)$ and the corresponding ROM solution $\boldsymbol{\tilde u}(\mu)$, the error is computed as
\begin{equation}
    \label{eq:relative-error}
    e(\mu) = \frac{\lVert \boldsymbol{u}(\mu) - \boldsymbol{\tilde u}(\mu) \rVert_{l^2}}{\lVert \boldsymbol{u}(\mu)  \rVert_{l^2}}.
\end{equation}
The relative errors for the two test cases are reported in \cref{tab:backstep-error,tab:cavity-error}, respectively, for both train and test configurations, in terms of average, standard deviation and median. \cref{fig:error-fields} shows the FOM solution for an unseen configuration, and the error fields for the baseline and novel methods. In particular, we only include the \quadnetmu error result, since, from \cref{tab:backstep-error,tab:cavity-error} it is the most accurate method.

\begin{table}[!htb]
\centering
    \caption{Comparison of models performance on the backstep dataset.}
    \begin{tabular}{lllllll}\toprule
    &\multicolumn{3}{c}{Train error}&\multicolumn{3}{c}{Test error}\\
    \cmidrule(lr){2-4} \cmidrule(lr){5-7}
    & \multicolumn{1}{c}{$\bar e$} & \multicolumn{1}{c}{$\sigma_e$} & \multicolumn{1}{c}{$\operatorname{Median}(e)$} & 
    \multicolumn{1}{c}{$\bar e$} & \multicolumn{1}{c}{$\sigma_e$} & \multicolumn{1}{c}{$\operatorname{Median}(e)$}\\ \midrule
    \podrbf         & 0.015       & 0.023       & 0.0059         & 0.015       & 0.022       & 0.006      \\
    \quadls         & 0.005       & 0.013       & 0.0015         & 0.0033      & 0.0062      & 0.0015     \\
    \quadnet        & 0.006       & 0.012       & 0.0015         & 0.0060       & 0.0099      & 0.0015     \\
    \quadnetmu      & \textbf{0.0014}      & 0.0029      & 0.00036        & \textbf{0.0013}      & 0.0023      & 0.00036    \\
    \bottomrule
    \end{tabular}
    \label{tab:backstep-error}
\end{table}

\begin{table}[!htb]
\centering
    \caption{Comparison of models performance on the cavity dataset.}
    \begin{tabular}{lllllll}\toprule
    &\multicolumn{3}{c}{Train error}&\multicolumn{3}{c}{Test error}\\
    \cmidrule(lr){2-4} \cmidrule(lr){5-7}
    & \multicolumn{1}{c}{$\bar e$} & \multicolumn{1}{c}{$\sigma_e$} & \multicolumn{1}{c}{$\operatorname{Median}(e)$} & 
    \multicolumn{1}{c}{$\bar e$} & \multicolumn{1}{c}{$\sigma_e$} & \multicolumn{1}{c}{$\operatorname{Median}(e)$}\\ \midrule
    \podrbf          & 0.030       & 0.010       & 0.027      & 0.033       & 0.010        & 0.030       \\
    \quadls          & 0.0154      & 0.0055      & 0.014      & 0.0176      & 0.0052      & 0.015      \\
    \quadnet         & 0.0154      & 0.0062      & 0.013      & 0.0182      & 0.0057      & 0.017      \\
    \quadnetmu       & \textbf{0.0091}      & 0.0035      & 0.0086     & \textbf{0.0108}      & 0.0032      & 0.010       \\
    \bottomrule
    \end{tabular}
    \label{tab:cavity-error}
\end{table}
The results show that in both test cases quadratic corrections improve the accuracy of the baseline linear ROM. The best model is \quadnetmu, which is more accurate than \quadls and \quadnet both on train and test data, and in both test cases. This is not surprising, as \quadnetmu is local in parameter space, whereas the other models are global. The advantage of a local model is especially evident for the backstep test case, where both the mean error and the standard deviations are substantially smaller. In fact, for this case the flow field shows rather different behaviours in the range of parameters considered, and the smaller variance of \quadnetmu indicates that this model is able to capture these behaviours better than the other two models. This, by contrast, is not so clear in the case of the lid-driven cavity, where the flow is \RA{convection dominated} over the whole range of parameters, and thus, although it is less accurately reproduced with a linear manifold, the flow structure does not change as much. 

From the tables we can also see that \quadnet is not necessarily better than \quadls, which, again, is something we expected.

Moreover, this space-continuous model can be trained with only a reduced number of selected sample points, making the training more efficient while maintaining a good accuracy, as we will see in \cref{subsec:partial}.

\subsection{Partial data}
\label{subsec:partial}
In this section we explore the possibility of exploiting the space-continuity of NN-based models to reduce the number of collocation points during training, and possibly reduce the computational cost of the offline phase. Given that \quadnetmu has shown the best performance, in what follows we will employ this model.

Training on a subset of collocation points introduces a further layer of complexity in the model, since now the question arises of how to choose the best subset of points given the initial mesh. The simplest strategy would be to randomly select points from a uniform distribution, so as to cover most of the domain. Although this strategy certainly works, an inspection of the correction terms reveals that often the corrections are quite localised, in the sense that they are mostly zero except for small regions of the domain. This suggests that choosing points uniformly in the domain may be sub-optimal, since we would rather have more points where the correction is large and changing rapidly, and fewer points where it is flatter and almost zero.

This result may be achieved by assigning each mesh node a Boltzmann-like probability of being chosen. More precisely, the probability of choosing the $i$-th node can be written as
\RA{
\begin{equation}
    p_i \propto \exp(-\tilde \tau /(\bar \tau_i + \varepsilon)),
\end{equation}
where $\bar \tau_i = \frac{1}{N_\mu} \sum_{j=1}^{N_\mu} \tau^{exact}(\boldsymbol x_i, \mu_j)$ is the average value of the correction at node $i$ across all training snapshots, and $\tilde \tau$ }is a suitable normalisation (for instance the average norm of the corrections, or the maximum  correction), and $\varepsilon$ is a small constant to avoid division by zero. \RA{This computation is easily carried out offline right after the exact corrections $\tau^{exact}$ are computed as in \cref{eq:corr-exact}}.

\begin{figure}[!htb]
    \centering
    \includegraphics[width=\linewidth,trim={4.5cm 0 0cm 0}, clip]{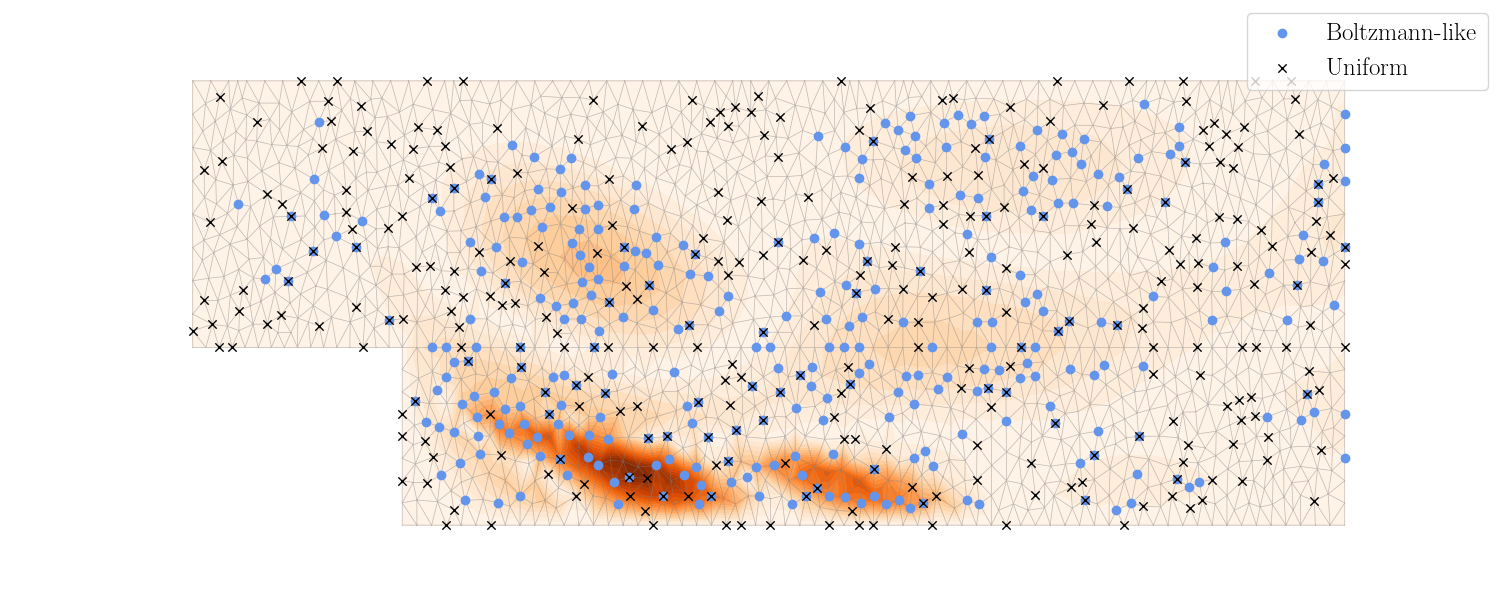}
    \caption{Comparison between samples of points drawn from different distribution. The blue points are sampled from a Boltzmann-like distribution, while the black crosses are sampled from a uniform distribution. A correction field is plotted in the background for reference.}
    \label{fig:sample}
\end{figure}

With such a distribution, points where the correction is larger are more likely to be chosen than points where the correction is typically smaller. The two sampling strategies are compared in \cref{fig:sample}, where we can see that the Boltzmann-like sampling does in fact yield the desired result.

\begin{figure}
    \centering
    \begin{subfigure}{0.48\textwidth}
        \includegraphics[width=\textwidth]{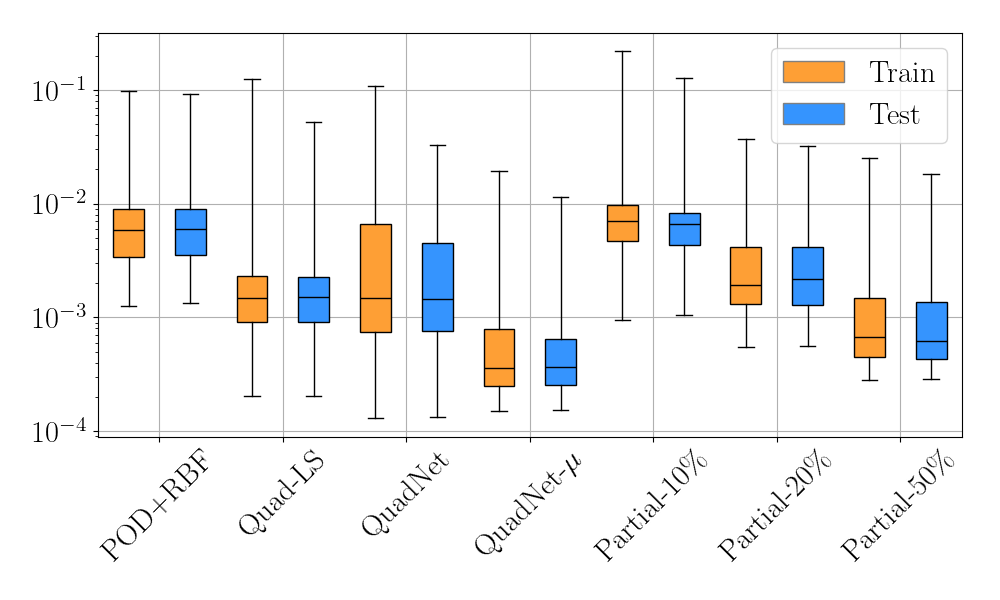}
        \caption{Comparison of errors for the backward-facing step.}
    \end{subfigure}
     \begin{subfigure}{0.48\textwidth}
        \includegraphics[width=\textwidth]{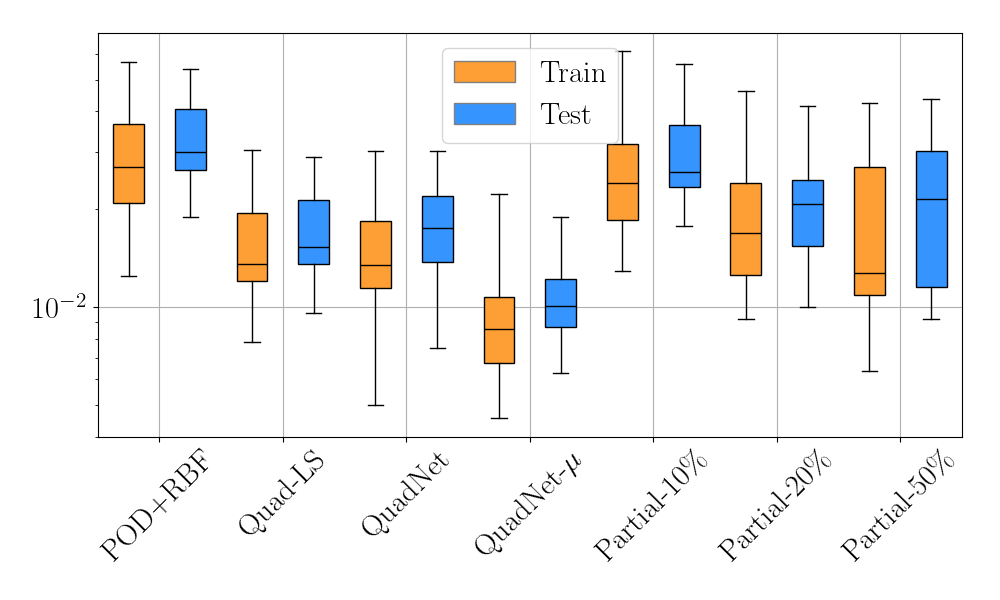}
        \caption{Comparison of errors for the lid-driven cavity.}
    \end{subfigure}
    \caption{Boxplots comparing the accuracy of models. \emph{Partial data} results refer to \quadnetmu. }
    \label{fig:boxplots-compare}
\end{figure}

We assess the applicability of \quadnetmu in the \emph{partial data} regime by training three models \RA{with the same architecture} on 10, 20 and 50\% of the mesh nodes, and comparing their performance to that of all models presented above in terms of the relative error \cref{eq:relative-error}. Results for both test cases are shown in \cref{fig:boxplots-compare}. We observe that the \emph{partial data} performance is comparable to that of the other models. As expected, there is a clear trend of decreasing error for increasing number of collocation points, with the best results achieved on the whole mesh and the worst on 10\% of mesh nodes. Moreover, for the lid-driven cavity the difference between the three \emph{partial data} cases is not as large as for the backward-facing step. This could be related to the complexity of the flow in this case, for which the field exhibits smaller features that are harder to capture with few collocation points. Finally, we believe it could be possible to improve the accuracy by sampling points with some other strategy, however this goes beyond the scope of this work and has not been explored further yet.

As for the goal of speeding up the offline phase, experiments have shown that training with fewer collocation points brings a reduction of roughly 10\% in the time required to train the networks. Since the problems considered here are fairly small, we expect an even larger reduction for more complex problems with larger meshes. Moreover, we noticed no significant differences in training times when using 10, 20 or 50\% of the points, which is something we can expect for this specific experimental setup. Indeed, on the one hand the problems here are so small that evaluating the network on 10, 20 or 50\% of the points makes very little difference when done in parallel. On the other hand, all tests were run on an Nvidia Quadro RTX 4000 GPU on a non-dedicated workstation, hence the time measurements may have been affected by unpredictable factors.

\subsection{Scarce data}
\label{subsec:scarce}
\begin{figure}[!htb]
    \centering
    \begin{subfigure}{0.8\textwidth}
        \centering
        \includegraphics[width=\linewidth,trim={4cm 0 1cm 0}, clip]{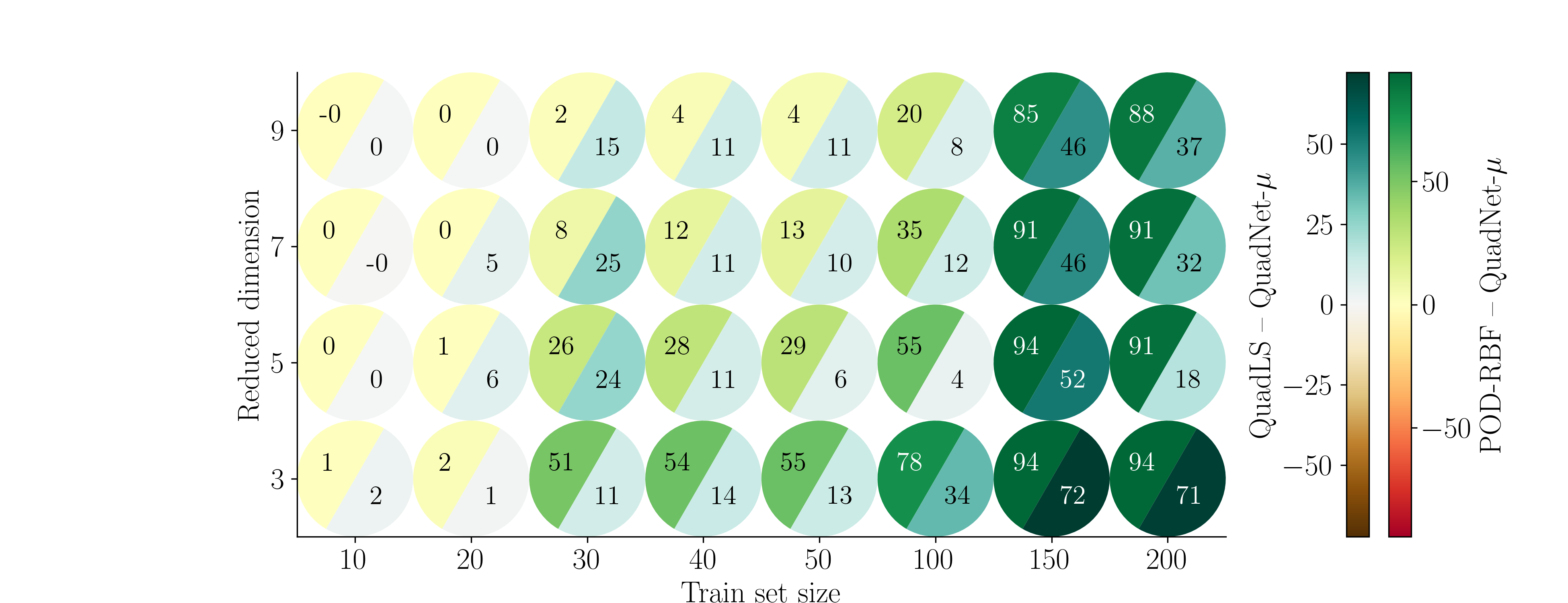}
        \caption{Backward-facing step.}
        \label{fig:backstep-heatmap}
    \end{subfigure}
    \begin{subfigure}{0.8\textwidth}
        \centering
        \includegraphics[width=\linewidth,trim={4cm 0 0cm 0}, clip]{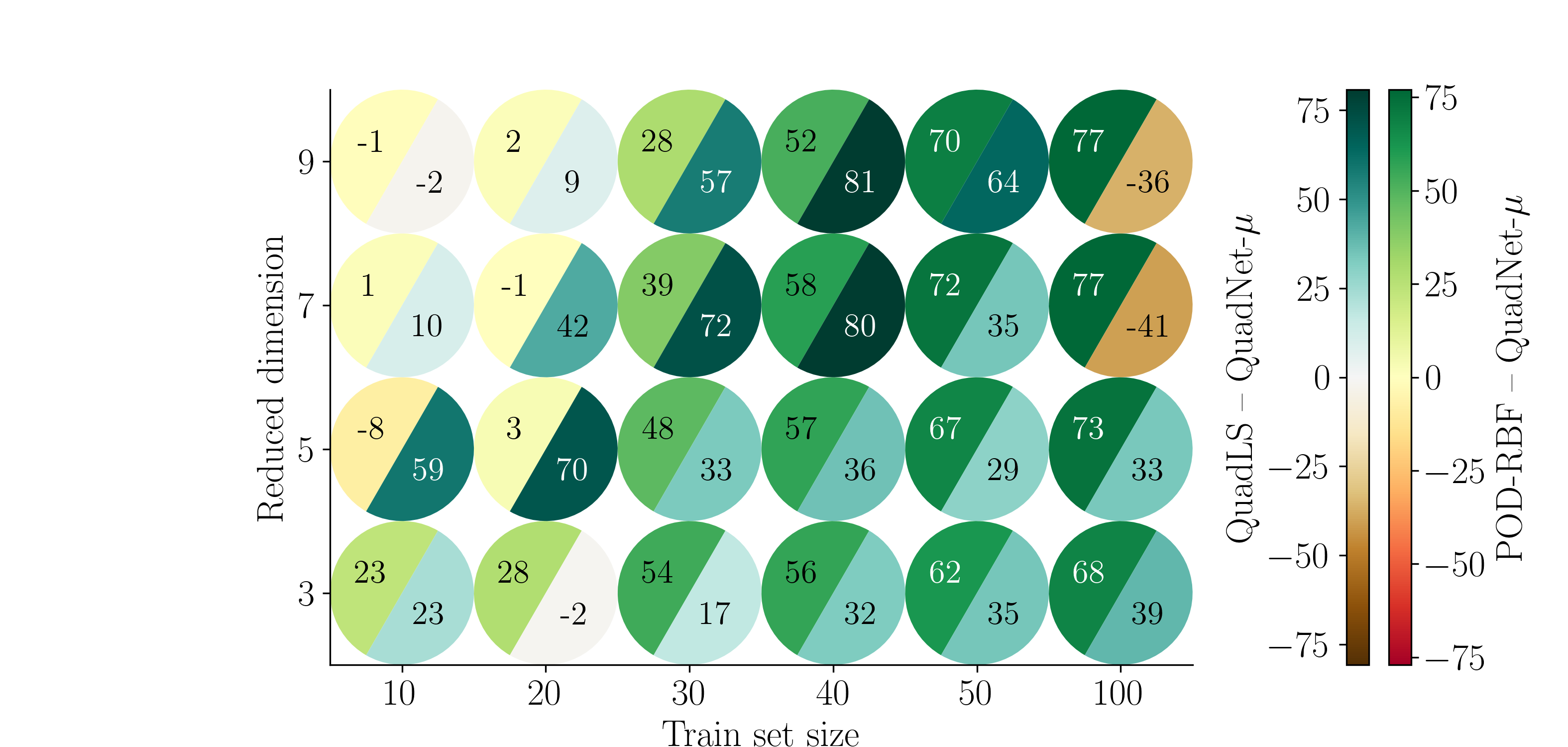}
        \caption{Lid-driven cavity.}
        \label{fig:cavity-heatmap}
    \end{subfigure}
    \caption{\RB{Comparison of \quadnetmu errors with \podrbf and \quadls. The left half of the circles shows the relative difference in error between \quadnetmu and \podrbf, while the right half shows the difference between \quadnetmu and \quadls.}}
    \label{fig:heatmaps}
\end{figure}

As a last test, we are interested in testing our model's performance when the training data are scarce. This scenario is especially common in industrial applications where the FOMs are too expensive to evaluate hundreds of times. In particular, we want to test \quadnetmu for different POD space dimensions and different training set sizes. As for the number of POD modes, we chose $r\in \{3,5,7,9\}$, while for the number of snapshots we chose $N_\mu \in \{10,20,30,40,50,100,150,200\}$ for the \RA{backward-facing} step and $N_\mu \in \{10,20,30,40,50,100\}$ for the lid-driven cavity. 
\RB{We trained \podrbf, \quadls and \quadnetmu for every combination of parameters and compared their errors on a test set.}
In \cref{fig:heatmaps} we show the \RB{relative error difference between \quadnetmu and both \podrbf and \quadls. Specifically, we compute}
\begin{equation}
    e_r = \frac{e_{\text{base}} - e_{\quadnetmu}}{e_{\text{base}}},
\end{equation}
\RB{where the subscript "base" is either "\podrbf" or "\quadls"}, and $e$ is computed as in \cref{eq:relative-error}.

We find that \quadnetmu performs better than the base model in almost all cases, with improvements in the error even larger than 90\%. We must note, however, that there are few cases where the quadratic model performs equally or even worse than the base model. These cases mostly correspond to those regimes where the number of snapshots is small and the number of POD modes is large. Indeed, in such cases there are two main issues. The first issue is that the POD representation of the solution manifold can be inaccurate if the few training snapshots are not representative of the whole manifold. In such cases, the error in the representation can be amplified by the quadratic model. The second issue, instead is related to the training of the model. Indeed, as $r$ approaches $N_\mu$, the POD representation of the training snapshots becomes more and more precise, which means that the corresponding correction term becomes smaller and smaller. As a consequence, the network learns to reproduce very small fields and thus the quadratic part of the corrected ROM is negligible in comparison to the linear part. In such cases we would expect the quadratic model to perform as good as the linear model.

\RB{As for the comparison with \quadls, we also find that \quadnetmu is generally more accurate, with the exception of very few cases most notably, for large datasets and relatively many POD modes. In general, we noticed that, since the networks grow in size as $r$ increases, their training becomes progressively harder, especially with few data available. In such cases we noticed a strong dependence on the choice of hyperparameters and, to some extent, on the weights initialisation. In this regard, special care should be taken in choosing an appropriate learning rate, possibly in combination with a learning rate scheduler. We believe that a fine tuning of the hyperparameters might lead to more accurate models. Moreover, we also point out that in cases where \quadnetmu is less accurate than \podrbf, it is still more accurate than \quadls. This means that in these cases, \quadls performs way worse than \podrbf. This can be clearly seen in \cref{fig:sensitivity_lines}, which shows a direct comparison of the errors from the three models.

As expected, the accuracy of all models is in general higher for larger training datasets and more POD modes. However, from the plots it is clear that \quadls struggles in providing meaningful corrections when the data are too scarce, thus yielding prediction errors that are substantially larger than those of both \podrbf and \quadnetmu. As the number of training snapshots increases, both \quadls and \quadnetmu improve the prediction accuracy, and for most combinations of parameters \quadnetmu outperforms \quadls. Moreover, we point out that \quadnetmu seems to be more robust than \quadls, in that the errors seem to be less spread out, as can be seen from the width of the shaded areas in \cref{fig:sensitivity_lines}. Notice that in all plots we show the median error together with the 5th and 95th percentiles. Similar considerations can be made by comparing the mean error instead, although for all models the presence of outliers makes the variance too large to provide a significant comparison. Moreover, we should note that in few cases, like $r=7$ or $r=9$ with large training size, \quadls slightly outperforms \quadnetmu, but always keeping very similar relative errors.

\begin{figure}[!htb]
    \centering
    \begin{subfigure}{.7\textwidth}
        \centering
        \includegraphics[width=\linewidth,trim={0cm 0 0cm 0}, clip]{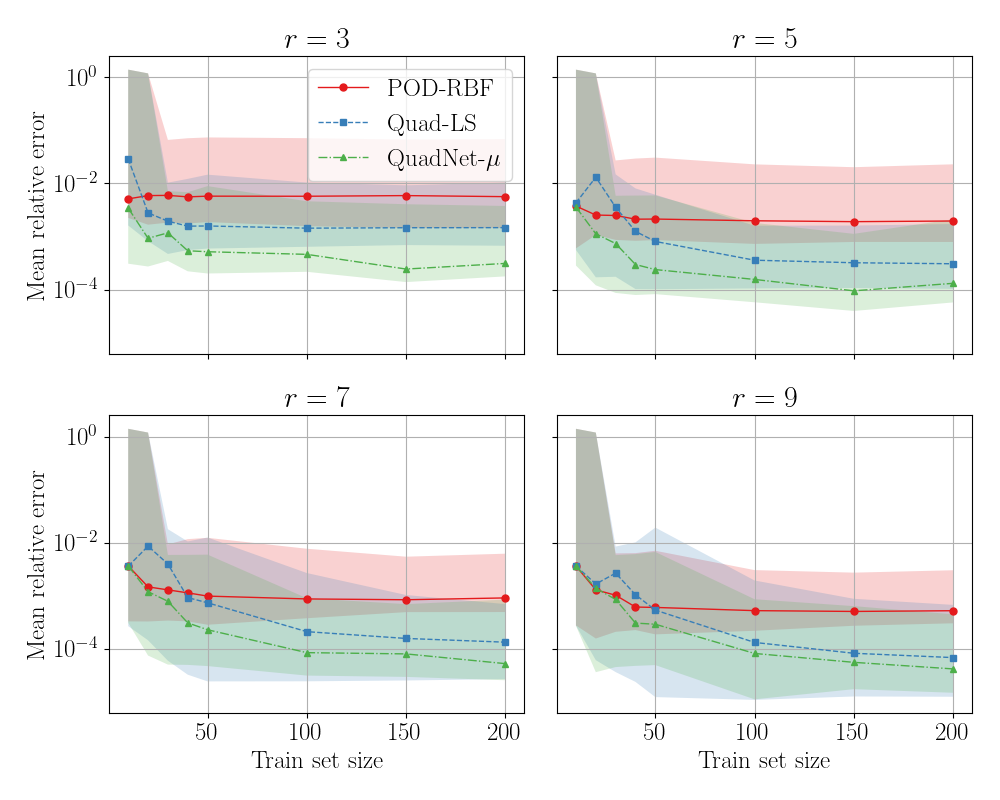}
        \caption{Backward-facing step.}
        \label{fig:backstep-lines}
    \end{subfigure}
    \begin{subfigure}{.7\textwidth}
        \centering
        \includegraphics[width=\linewidth,trim={0cm 0 0cm 0}, clip]{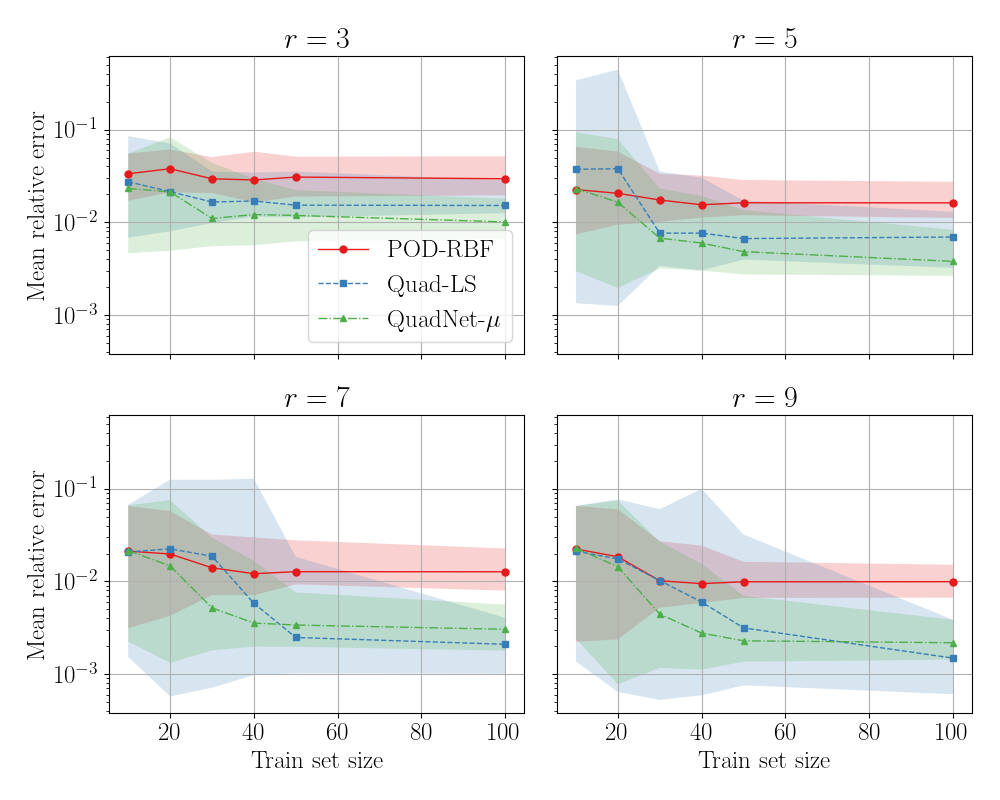}
        \caption{Lid-driven cavity.}
        \label{fig:cavity-lines}
    \end{subfigure}
    \caption{Comparison between \podrbf, \quadls and \quadnetmu errors for different combinations of $r$ and $N_\mu$. The markers refer to the median relative error, while the shaded area is between the 5th and 95th percentiles.}
    \label{fig:sensitivity_lines}
\end{figure}
}

\section{Conclusions}\label{sec:conclusions}

In this paper, we presented a novel class of non-intrusive reduced-order models (ROMs) for computational fluid dynamics that employs deep neural operator networks to overcome the limitations inherent to classical linear and quadratic closure models. Specifically, we introduced two neural-network-based methods, QuadNet and QuadNet-$\mu$, built upon the Deep Operator Network (DeepONet) and Multi-Input Operator Network (MIONet), respectively, to learn a quadratic correction for a given surrogate model.

The proposed methods demonstrated substantial improvements in accuracy over traditional linear POD-based methods and quadratic interpolation schemes such as POD-RBF. In particular, numerical experiments on benchmark fluid dynamics test cases confirmed that our methods achieve up to 90\% error reduction compared to POD-RBF models\RB{, as well as up to 80\% error reduction compared to \quadls models}, while requiring fewer training data points and fewer model parameters. These results underline the effectiveness of neural operators in capturing complex nonlinear interactions within fluid flows.
The integration of continuous mappings in both spatial coordinates and parameter space further enhanced the flexibility and generalizability of the proposed models, making them suitable for practical scenarios characterized by limited data availability. 

Future work will explore extensions of our methodology to even more complex problems, investigating the extension of this framework to higher-order correction. The neural architecture can indeed be easily changed to learn operators of higher order, which may lead to more accurate solutions.
Additionally, further studies will assess the robustness and efficacy of the proposed framework when combined with more sophisticated base models such as POD-NN, aiming to consistently maintain the demonstrated advantages in terms of accuracy and generalization.

\bibliographystyle{plain}
\bibliography{refs}

\end{document}